\documentclass[article,12pt]{article}
\usepackage{graphicx}        % standard LaTeX graphics tool
\usepackage{amsmath,amssymb,amsfonts}        % AMS math packages
\newcommand{\ud}{\mbox{d}}
\usepackage{bm}
\newcommand{\iprod}[1]{\langle#1\rangle}
\usepackage{color}
\usepackage{framed}
\usepackage{lipsum}
\renewenvironment{shaded}{%
  \MakeFramed{\advance\hsize-\width \FrameRestore\FrameRestore}}%
 {\endMakeFramed}
\definecolor{shadecolor}{gray}{0.75}

\begin{document}
\title{Wider contours and adaptive contours}
\author{Shev MacNamara\footnote{School of Mathematical and Physical Sciences, University of Technology, Sydney, \texttt{shev.macnamara@uts.edu.au}}, \, William McLean\footnote{The School of Mathematics and Statistics, University of New South Wales, Sydney, \texttt{w.mclean@unsw.edu.au}} \,  and   Kevin Burrage\footnote{ ARC Centre of Excellence for Mathematical and Statistical Frontiers, School of Mathematical Sciences, Queensland University of Technology, Brisbane, Australia
and
Visiting Professor to the University of Oxford, U.K.,
\texttt{kevin.burrage@qut.edu.au}}}
\maketitle

\abstract{Contour integrals in the complex plane are the basis of effective numerical methods for computing matrix functions, such as the matrix exponential and the Mittag-Leffler function.
These methods provide successful ways to solve partial differential equations, such as convection--diffusion models.
Part of the success of these methods comes from exploiting the freedom to choose the contour, by appealing to Cauchy's theorem.   
However, the pseudospectra of non-normal matrices or operators present a challenge for these methods: if the contour is too close to regions where the norm of the resolvent matrix is large, then the accuracy suffers.
Important applications that involve non-normal matrices or operators include the Black--Scholes equation of finance, and Fokker--Planck equations for stochastic models arising in biology.
Consequently, it is crucial to choose the contour carefully. 
As a remedy, we discuss choosing a contour that is wider than it might otherwise have been for a normal matrix or operator.
We also suggest a semi-analytic approach to adapting the contour, in the form of a parabolic bound that is derived by estimating the field of values.   
To demonstrate the utility of the approaches that we advocate, we study three models in biology: a  monomolecular reaction, a  bimolecular reaction and a  trimolecular reaction.
Modelling and simulation of these reactions is done within the framework of Markov processes. 
We also consider non-Markov generalisations that have Mittag-Leffler waiting times instead of the usual exponential waiting times of a Markov process.
}

\section{Introduction}
\label{sec:1}

We begin with the Chapman--Kolmogorov forward equation, associated with a Markov process on discrete states, for the evolution in continuous time of the probability of being in state $j$ at time $t$:
\begin{eqnarray}
 \frac{\ud }{\ud t}  p(j,t)  &=&   -  |a_{jj}|p(j,t)   + \sum_{i \ne j} a_{j,i} p(i,t). 
\label{eq:markov:chapman-kolmogorov}
\end{eqnarray}
Here for $j \ne i$, $a_{i,j} \ge 0$, and $a_{i,j} \ud t$ is approximately the probability to transition from state $j$ to state $i$ in a small time $\ud t$.
The diagonal entry  of  an associated matrix $\mathbb{A} = \{ a_{i,j}\}$ is  defined by the requirement that the matrix has columns that sum to zero, namely $a_{jj} = - \sum_{i \ne j} a_{i,j}$.
This equation in the Markov setting \eqref{eq:markov:chapman-kolmogorov} can be generalised to a non-Markovian form:
\begin{eqnarray}
 & &\frac{\ud }{\ud t}  p(j,t)  = \nonumber \\
 & &   -  \int_0^{t} K(j,t-u) p(j,u)  \ud u + \sum_{i \ne j} \frac{a_{j,i}}{|a_{jj}|}  \int_0^{t} K(i,t-u) p(i,u)  \ud u. \nonumber \\
\label{eq:master:equation:general}
\end{eqnarray}
Here the so-called memory function $K(i,t-u)$ is defined via its Laplace transform $\hat{K}(j,s) $, as the ratio of the Laplace transform of the waiting time to the Laplace transform of the survival time.
The process is a Markov process if and only if the waiting times are exponential random variables, in which case the $K(j,t)$ appearing in the convolutions in \eqref{eq:master:equation:general} are Dirac delta distributions and \eqref{eq:master:equation:general} reduces to the usual equation in \eqref{eq:markov:chapman-kolmogorov}.
In the special case of Mittag-Leffler waiting times, \eqref{eq:master:equation:general} can be re-written as \cite{ShevCauchyIntegralMasterEqnPseudoSpectraCTAC2015,MacNamaraFractionalEulerLimit2016}:
\begin{equation}
D_t^{\alpha} \bm{p} = \mathbb{A} \bm{p}  \qquad \quad \textrm{\;\; with solution} \qquad \; \bm{p}(t) = E_{\alpha}(\mathbb{A} t^{\alpha}) \bm{p}(0).
\label{eq:mittag:leffler:master:equation}
\end{equation}
Here $D_t^{\alpha} $ denotes the Caputo fractional derivative, and the solution comes in terms of a Mittag-Leffler function 
\begin{equation}
 E_{\alpha}(z) = \sum_{k=0}^{\infty} \frac{z^{k}}{\Gamma(\alpha k +1)} .
 \label{eq:mittag:leffler}
\end{equation}
Here $0 \le \alpha \le 1$.
In Section \ref{sec:ML:ssa}, we provide a \texttt{MATLAB} code to simulate sample paths of this stochastic process with Mittag-Leffler waiting times.
When $\alpha=1$ the series \eqref{eq:mittag:leffler} reduces to the usual exponential function, and the fractional equation in \eqref{eq:mittag:leffler:master:equation} reduces to the original  Markov process \eqref{eq:markov:chapman-kolmogorov}.
In Section \ref{sec:ML:contour:integral}, we provide a \texttt{MATLAB} code to compute the solution $  E_{\alpha}(\mathbb{A} t^{\alpha}) \bm{p}(0)$ in \eqref{eq:mittag:leffler:master:equation} directly, via a contour integral.

Next, we introduce the (Markovian) Fokker--Planck partial differential equation (sometimes also known as the forward Kolmogorov equation)
\begin{equation}
\frac{\partial }{ \partial t} p(x,t) = - \frac{\partial }{ \partial x} \left( a(x) p(x,t)  \right) + \frac{1}{2}\frac{\partial^2 }{ \partial x^2} \left( b(x) p(x,t)  \right)
\label{eq:fokker-planck}
\end{equation}
for a probability density $p(x,t)$ in one space dimension and time, with suitable boundary conditions, initial conditions, and smoothness assumptions on the coefficients $a(x)$, and on $b(x) \ge 0$.
Later, we will use complex-variable methods, so it is worth noting at the outset that our coefficients  $a(x)$ and $b(x)$ are always real-valued.
It is also worth noting that $a(x)$ here in \eqref{eq:fokker-planck} is not the same as $a_{i,j}$ appearing in the matrix above in \eqref{eq:mittag:leffler:master:equation}, although there is a close relationship that allows one to be deduced from the other.
This PDE for the density corresponds to a stochastic differential equation for sample paths 
\begin{equation}
\ud X =  a(X) \ud t +  \sqrt{b(X)} \ud W.
\label{eq:SDE}
\end{equation}
Both the PDE \eqref{eq:fokker-planck} and SDE \eqref{eq:SDE} are continuous in the spatial variable.
An introduction to these models, and their connections with discrete versions, can be found in \cite{AndersonKurtzChapter2011,GillespieIsomerization2002}.
Our discrete models do respect important properties such as non-negativity.
However, there are issues with the Fokker--Planck PDE model, the Chemical Langevin SDE model, and other diffusion approximations: often these approximations do not maintain positivity.
 These issues are discussed by Leite \& Williams  in the Kolmogorov Lecture and accompanying paper, where a method that maintains non-negativity is proposed \cite{WilliamsKolmogorovLecture2016}.

We do not simulate or solve either of the PDE \eqref{eq:fokker-planck} or the SDE \eqref{eq:SDE}.
We do however simulate and solve the closely related models that are discrete in space and that are governed by master equations \eqref{eq:markov:chapman-kolmogorov} or generalized master equations \eqref{eq:master:equation:general}, which can be thought of as finite difference discretizations of  \eqref{eq:fokker-planck}.
In particular, the PDE \eqref{eq:fokker-planck} can be thought of as a continuous-in-space analogue of the discrete process in  \eqref{eq:markov:chapman-kolmogorov} that involves a matrix $\mathbb{A}$. 
The utility of the PDE is that it is easier to find an estimate of the field of values of the spatial differential operator on the right hand side of \eqref{eq:fokker-planck} than of the corresponding matrix $\mathbb{A}$. 
We can then use an estimate of one as an approximation for the other.

Developing appropriate methods for modelling and simulation of these important stochastic processes is a necessary first step for more advanced scientific endeavors.
A natural next step is an inverse problem, although we do not pursue that in this article. 
For example, it is of great interest to estimate the rate constants in the models described in Section~\ref{sec:3:Models}, typically based on limited observations of samples resembling the simulation displayed in Figure~\ref{fig:ML:SSA}. 
It is more challenging to estimate the parameter $\alpha$ in fractional models, or indeed to address the question of whether or not a single $\alpha$ (the only case considered in this article) is appropriate.
Interest in fractional models is growing fast as they are finding wide applications including in cardiac modelling \cite{BuenoKayGrauRodriguezBurrage2014}, and this includes exciting new applications of Mittag-Leffler functions \cite{BuenoBurrage2016}.

Section~\ref{sec:3:Models} introduces three models that can be cast in the form of \eqref{eq:mittag:leffler:master:equation} and Section~\ref{sec:pseudospectra} uses these models as vignettes to exhibit manifestations of pseudospectra.
Next we introduce two methods of simulation for these models.
A Monte Carlo approach is presented in Section~\ref{sec:ML:ssa}. 
Section~\ref{sec:ML:contour:integral} presents an alternative method that directly computes the solution of \eqref{eq:mittag:leffler:master:equation} as a probability vector  via a contour integral.
Finally we suggest a bound on the field of values that is useful when designing contour methods.

To our knowledge Figure~\ref{fig:schlogl:pseudospectra} is the first visualization of the pseudospectra of the Schl\"ogl reactions. 
The estimate in \eqref{eq:mono:molecular:matrix:parabolic:bound}  is also new.
In fact \eqref{eq:mono:molecular:matrix:parabolic:bound} is the specialization of our more general estimates appearing in~\eqref{def:zero:dirichlet:bound} and~\eqref{def:zero:flux:neumann:bound} to the monomolecular model, but it should be possible to likewise adapt our more general estimates to other models such as bimolecular models.

\section{Three fundamental models}
\label{sec:3:Models}
All three models that we present are represented by \textit{tri-diagonal matrices}: 
$
 j \notin \{ i-1,i, i+1\}  \Rightarrow  \mathbb{A}_{i,j}=0.
$
%Non-zeros entries appear on the main diagonal, and on the two immediately adjacent diagonals of the matrix.
Since all other entries are zero, below, we only specify the non-zero entries on the three main diagonals.
In fact,  an entry on the main diagonal is determined by the requirement that columns sum to zero (which  corresponds to conservation of probability), so it would suffice to specify only the two non-zero diagonals immediately below and above the main diagonal.

\subsection{Monomolecular, bimolecular and trimolecular models}

A model of a monomolecular reaction 
\[
S_1 \leftrightarrows S_2
\]
such as arises in chemical isomerisation \cite{GillespieIsomerization2002}, or in models of ion channels in cardiac-electrophysiology, or neuro-electrophysiology, can be represented by the following $N \times N$ \textbf{ monomolecular model} matrix~
\begin{equation}
 \mathbb{A}_{i,j} =
  \begin{cases}
    j-1, & i=j-1,\\
    -m, & i=j,\\		
    m-j+1, & i=j+1.
  \end{cases}
\label{eq:mono:molecular:matrix}
\end{equation}

Here we have assumed that the rate constants are equal to unity, $c_1=c_2=1$, and we have assumed that $m=N-1$ where $m$ is the maximum number of molecules.
More details, including exact solutions, can be found in  \cite{IserlesMacNamara2017}.
An instance of this matrix when $m=5$ is
  \begin{equation}
\mathbb{A} =
\left(
\begin{tabular}{rrrrrr}
   $-5$    & $1$ &    \\
   $5$   & $-5$ &    $2$ \\
         & $4$ &  $-5$ &    $3$  \\
       &     &   $3$  &  $-5$  &  $4$  &   \\
        &    &    &    $2$ &   $-5$  &   $5$ \\
        &     &     &     &  $1$  &  $-5$
   \end{tabular}
\right).
\label{eq:A}
\end{equation}
The model is two-dimensional, but the conservation law allows it to be treated as effectively one-dimensional, by letting $x$ represent $S_1$, say.
Then one possible corresponding continuous model  \eqref{eq:fokker-planck} has drift coefficient $a(x)=-c_1x+c_2(m-x)$ and diffusion coefficient  $b(x)=c_1x+c_2(m-x)$, for $0<x<m$.
In the discrete Markov case, the exact solution is binomial.
When $c_1=c_2$ the stationary probability density of the continuous Markov model is given by Gillespie  \cite{GillespieIsomerization2002} as a Gaussian, which highlights the issues associated with  the continuous models such as choosing the domain, boundary conditions, and respecting positivity. 
To enforce positivity, we might instead pose the PDE on the  restricted domain $0 \le x \le m$.

Next we introduce a model for the bimolecular reaction 
\[
S_1 + S_2 \leftrightarrows S_3,
\]
via the following $N \times N$ \textbf{ bimolecular model} matrix
\begin{equation}
 \mathbb{A}_{i,j} =
  \begin{cases}
    j-1, & i=j-1,\\
    -j+1 -(m-j+1)^2 & i=j,\\		
    (m-j+1)^2, & i=j+1.
  \end{cases}
\label{eq:bi:molecular:matrix}
\end{equation}
~\\

This matrix and the model are also introduced in \cite{ShevCauchyIntegralMasterEqnPseudoSpectraCTAC2015}, where more details can be found.
A small example with $m=5$ is 
  \begin{equation}
A =
\left(
\begin{tabular}{rrrrrr}
   $-25$    & $1$ &    \\
   $25$   & $-17$ &    $2$ \\
         & $16$ &  $-11$ &    $3$  \\
       &     &   $9$  &  $-7$  &  $4$  &   \\
        &    &    &    $4$ &   $-5$  &   $5$ \\
        &     &     &     &  $1$  &  $-5$
   \end{tabular}
\right).
\label{eq:A:bimolecular}
\end{equation}
Here we have assumed that the rate constants are equal to unity, $c_1=c_2=1$, and that the initial condition  is $[S_1,S_2,S_3] = [m,m,0]$, so $m$ is the maximum number of molecules, and $m=N-1$.
The model is three-dimensional, but the conservation law together with this initial condition allow it to be treated as effectively one-dimensional, by letting $x$ represent $S_3$, say.
Then a possible corresponding continuous model  \eqref{eq:fokker-planck} has drift coefficient $a(x)=-c_1x+c_2(m-x)^2$ and diffusion coefficient  $b(x)=c_1x+c_2(m-x)^2$.

Finally, we introduce the Schl\"ogl model \cite{MacNamaraFractionalEulerLimit2016}, which consists of two reversible reactions
\[
B_1 + 2X  \longleftrightarrow 3X, \qquad \qquad 
B_2 \longleftrightarrow X.
\]
Here $B_1=1 \times 10^{5}$, $B_2=2 \times 10^{5}$. 
%More details can be found in \cite{MacNamaraFractionalEulerLimit2016}.
%Note that a genuinely  trimolecular reaction mechanism is for physical reasons unlikely; it is more plausible that a  trimolecular scheme is a good approximation to two consecutive  bimolecular reactions that involve an unstable intermediate.
The associated matrix is given below, where $k_1 = 3 \times 10^{-7}$, $k_2 = 1 \times 10^{-4}$, $k_3 = 1 \times 10^{-3}$, and $k_4 = 3.5$.
In theory this matrix is infinite but we truncate to a finite section for numerical computation.

We define the Schl\"ogl model matrix (an example of a \textbf{ trimolecular model} scheme)~\\
\begin{equation}
 \mathbb{A}_{i,j} =
  \begin{cases}
    \frac{1}{6}k_2(j-1)(j-2)(j-3)+k_4(j-1), & i=j-1,\\		
    k_3B_2 + \frac{1}{2}k_1B_1(j-1)(j-2), & i=j+1.
  \end{cases}
\label{eq:tri:molecular:matrix}
\end{equation}
% If we index rows and columns starting at j==0, 
% then the matrix is as follows...
%\begin{equation}
% \mathbb{A}_{i,j} =
%  \begin{cases}
%    \frac{1}{6}k_2j(j-1)(j-2)+k_4j, & i=j-1,\\
%    -(  \frac{1}{6}k_2j(j-1)(j-2)+k_4j, +   k_3B_2 + \frac{1}{2}k_1B_1j(j-1)), & i=j,\\		
%    k_3B_2 + \frac{1}{2}k_1B_1j(j-1), & i=j+1.
%  \end{cases}
%\label{eq:tri:molecular:matrix}
%\end{equation}

For $i=j$, the diagonal entry is $-(  \frac{1}{6}k_2(j-1)(j-2)(j-3)+k_4(j-1) +   k_3B_2 + \frac{1}{2}k_1B_1(j-1)(j-2))$.
The first column is indexed by $j=1$ and corresponds to a state with $0=j-1$ molecules.
The corresponding continuous model  \eqref{eq:fokker-planck} has drift coefficient $a(x)=  k_3B_2 + \frac{1}{2}k_1B_1x(x-1) - \frac{1}{6}k_2x(x-1)(x-2)-k_4x $ and diffusion coefficient  $b(x)= k_3B_2 + \frac{1}{2}k_1B_1x(x-1) + \frac{1}{6}k_2x(x-1)(x-2) +k_4x $.

\section{The significance of pseudospectra for stochastic processes}
\label{sec:pseudospectra}
All three matrices introduced in the previous section exhibit interesting pseudospectra.
As an illustration of the way that the pseudospectra manifest themselves, we will now consider numerically computing eigenvalues of the three matrices.
This is a first demonstration that one must respect the pseudospectrum when crafting a numerical method. 
That issue will be important again when we use numerical methods based on contour integrals to compute Mittag-Leffler matrix functions.

\begin{figure}[t]
%\sidecaption[t]
\centering
\includegraphics[scale=.7]{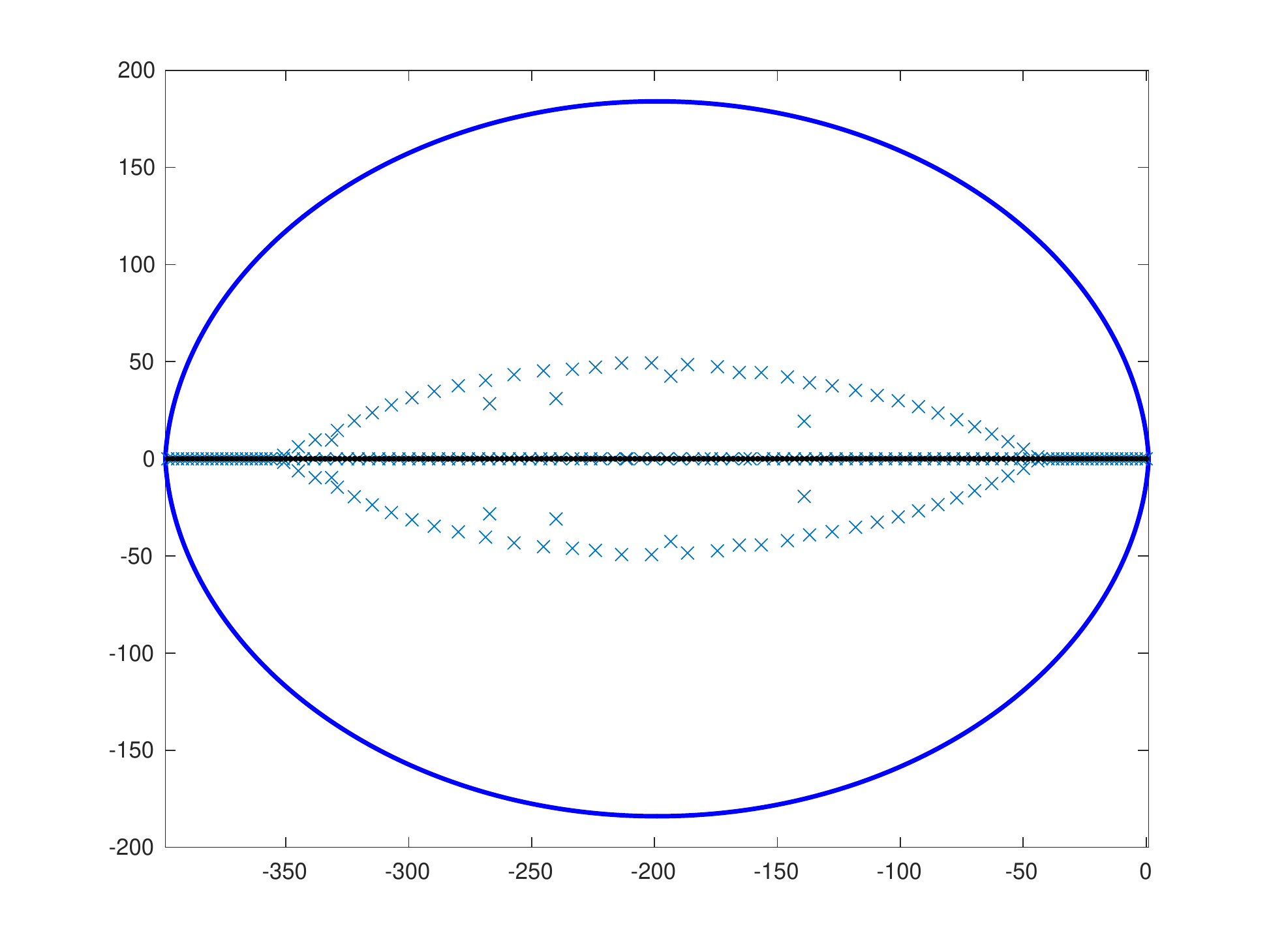}
\caption{The field of values (solid line resembling an oval-like shape) for the discrete  monomolecular matrix in \eqref{eq:mono:molecular:matrix} when $N=200$ ($m=N-1$), as computed by \texttt{Chebfun} \cite{Driscoll2014}.
The crosses mark numerically computed eigenvalues via \texttt{eig}, but complex eigenvalues are \textit{wrong.}
The true eigenvalues are purely real and are marked on the real axis by dots (note that the dots are so close together that they may seem to be almost a solid line).
These correct values can come by instead computing the eigenvalues of the symmetrized matrix, after using the diagonal scaling matrix created by the \texttt{MATLAB} listing provided here. 
}
\label{fig:MonoMolecular:FOV}       % Give a unique label
\end{figure}

\begin{figure}[t]
%\sidecaption[t]
\centering
\includegraphics[scale=.7]{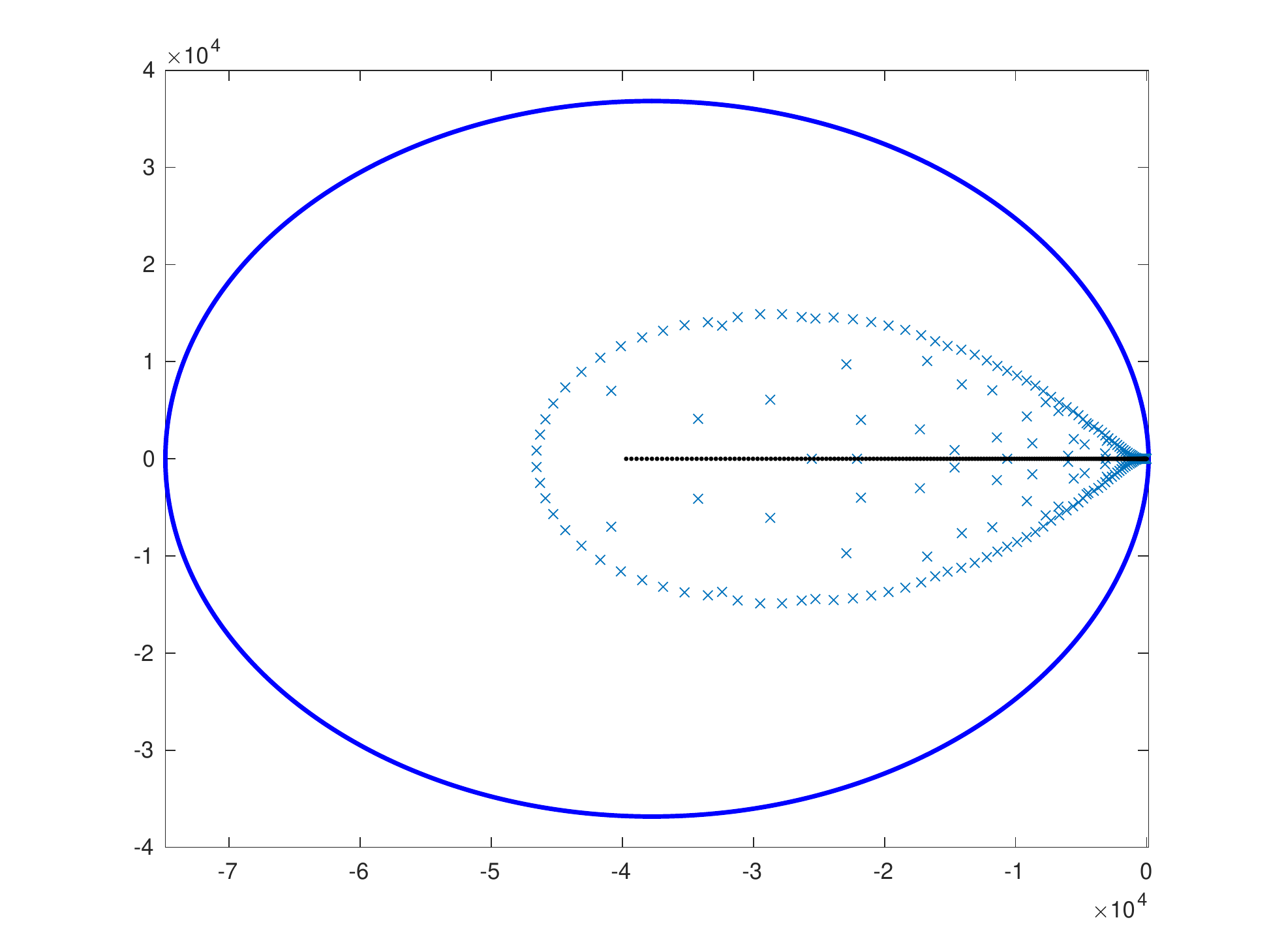}
\caption{Same as Figure~\ref{fig:MonoMolecular:FOV},  for the  bimolecular model.
The field of values  for the discrete  bimolecular matrix in \eqref{eq:bi:molecular:matrix} when $N=200$ is the solid line resembling an oval-like shape, as computed by \texttt{Chebfun} \cite{Driscoll2014}.
}
\label{fig:BiMolecular:FOV}       % Give a unique label
\end{figure}

\begin{figure}[t]
%\sidecaption[t]
\centering
\includegraphics[scale=0.9]{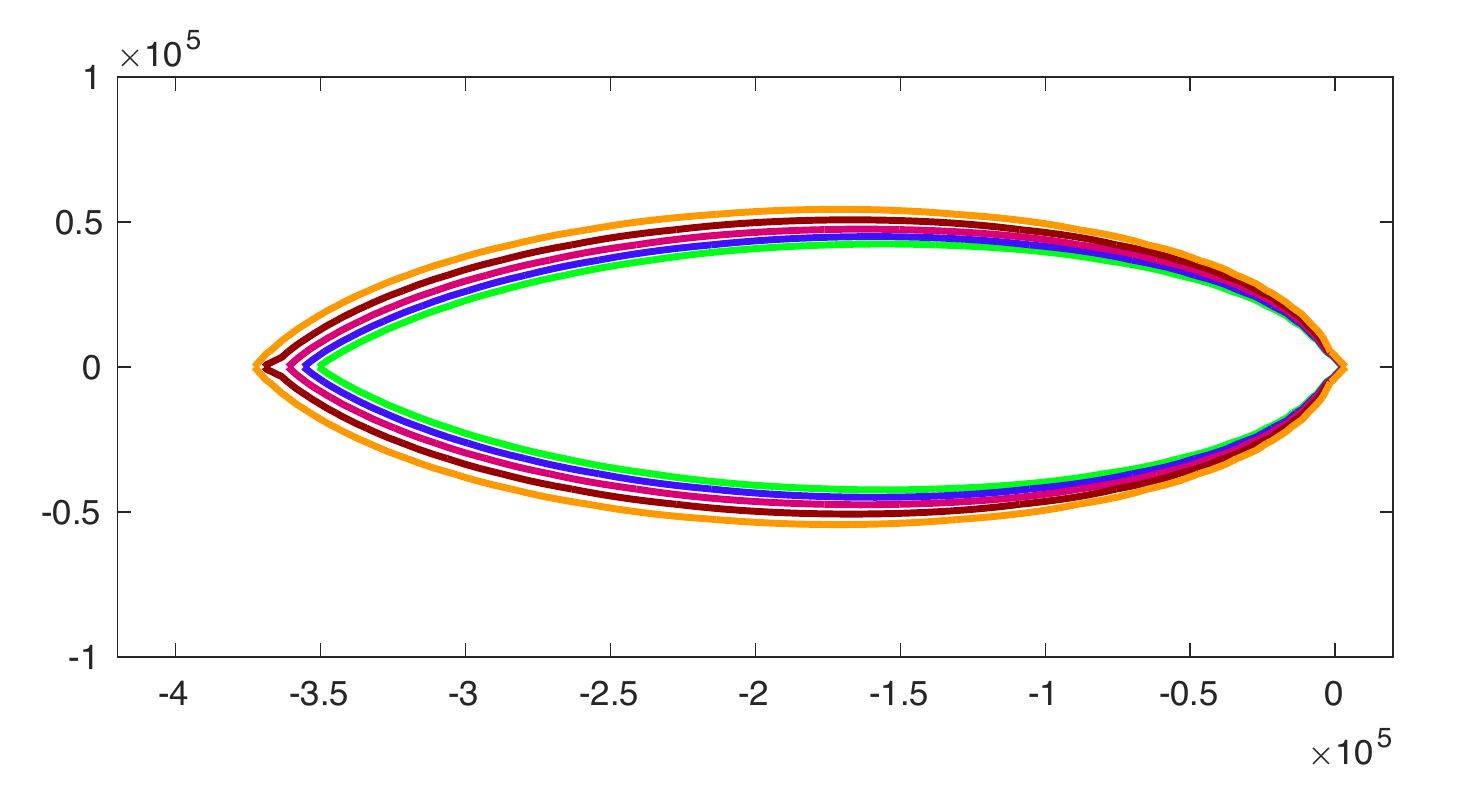}
\caption{Pseudospectrum of a $2000 \times 2000$ finite section of the matrix in \eqref{eq:tri:molecular:matrix} representing the Schl\"ogl reactions, as computed by  \texttt{EigTool}.
Contours correspond to $\epsilon=10^{-2} , 10^{-4}, 10^{-6}, 10^{-8}, 10^{-10}$ in \eqref{eq:pseudospectra:definition}.
The contour closest to the real axis corresponds to $10^{-10}$.
}
\label{fig:schlogl:pseudospectra}       % Give a unique label
\end{figure}

\begin{figure}[t]
%\sidecaption[t]
\centering
\includegraphics[scale=.7]{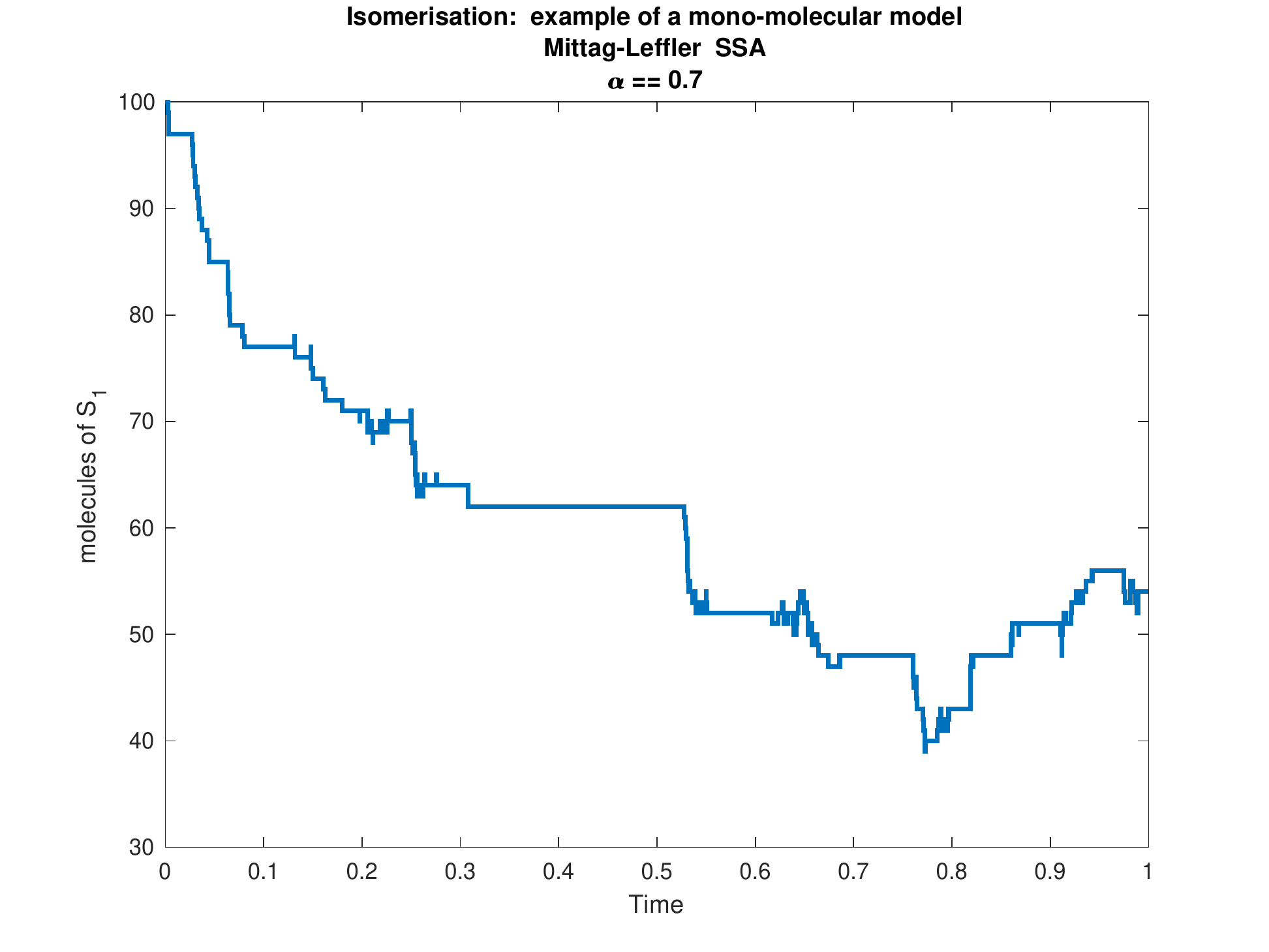}
\caption{A  sample simulation of the monomolecular reaction, with the Mittag-Leffler SSA provided here in the \texttt{MATLAB} listing.
Parameters: $\alpha = 0.7$, $t=1$, and initial condition a Dirac delta distribution on the state $[S_1,S_2]=[m,0]=[100,0]$. 
Compare with Figure~\ref{fig:ML:ME}.}
\label{fig:ML:SSA}       % Give a unique label
\end{figure}

The reader can readily verify that using any standard eigenvalue solver, such as \texttt{eig} in \texttt{MATLAB}, leads to numerically computed eigenvalues that are complex numbers.
However, these numerically computed complex eigenvalues are \textit{wrong:} the \textit{true eigenvalues are purely real.} 
It is the same story for all three models. 
See the numerically computed eigenvalues displayed here in Figure~\ref{fig:MonoMolecular:FOV} and Figure~\ref{fig:BiMolecular:FOV}, for example.
We suggest this effect happens much more widely for models arising in computational biology.

%As an example, consider the  monomolecular matrix \eqref{eq:mono:molecular:matrix}.
%Numerically computed eigenvalues for this matrix are displayed in Figure~\ref{fig:MonoMolecular:FOV}.
Figure~\ref{fig:MonoMolecular:FOV}  and  Figure~\ref{fig:BiMolecular:FOV} make use of the following method to create a diagonal scaling matrix.
Here is a \texttt{MATLAB} listing to create a diagonal matrix $D$ that symmetrizes a tridiagonal matrix of dimension $N$ of the form described in any of the three models considered here:

\begin{verbatim}
d(1) = 1;
for i = 1:N-1
    d(i+1) = sqrt(A(i,i+1)/A(i+1,i)) *d(i);
end
D = diag(d);   Asym = D*A*inv(D);
\end{verbatim}

This symmetrization by a diagonal matrix in a similarity transform is known to physicists as a gauge transformation, and it is described by Trefethen and Embree \cite[Section 12]{TreEmb05}.
A real symmetric matrix has real eigenvalues so the eigenvalues of $D A D^{-1}$ are purely real.
The matrix $D A D^{-1}$ and the matrix $A$ share the same eigenvalues because they are similar.
This is one way to confirm that the true eigenvalues of $A$ are purely real.
Numerical algorithms typically perform well on real symmetric matrices, so it is better to numerically compute the eigenvalues of $D A D^{-1}$ than to compute eigenvalues of $A$ directly.

The $\epsilon$-pseudospectra \cite{TreEmb05} can be defined as the region of the complex plane where the norm of the resolvent matrix is large: the set $z \in \mathbb{C}$ such that
\begin{equation}
||(z \mathbb{I} - \mathbb{A})^{-1} ||  > \frac{1}{\epsilon}.
\label{eq:pseudospectra:definition}
\end{equation}
We use the $2$-norm in this article.
Equivalently, this is the region of the complex plane for which $z$ is an eigenvalue of $(\mathbb{A} + \mathbb{E})$ for some  perturbing matrix that has a small norm: $||\mathbb{E} ||<\epsilon$.
Thus the wrong complex `eigenvalues'  appearing as crosses in  Figure~\ref{fig:MonoMolecular:FOV} and Figure~\ref{fig:BiMolecular:FOV}  offer a way to (crudely) visualise the pseudospectrum.
Numerical visualizations of the pseudospectra of the family of  monomolecular matrices defined by \eqref{eq:mono:molecular:matrix} can be found in \cite{IserlesMacNamara2017},  and the pseudospectra of the bimolecular reaction defined by \eqref{eq:bi:molecular:matrix} can be found in \cite{ShevCauchyIntegralMasterEqnPseudoSpectraCTAC2015}.
Here, as an example of a  trimolecular scheme, we present in Figure~\ref{fig:schlogl:pseudospectra}   the pseudopspectra for the Schl\"ogl reactions \eqref{eq:tri:molecular:matrix}, as computed by \texttt{EigTool}.
The resolvent norm is largest near the negative real axis, as we would expect because that is where the true eigenvalues are located.
The level curves are not equally spaced in the figure, becoming bunched up together, suggesting a more interesting structure that remains to be elucidated as the dimension of the matrix grows.
%Thus in the `eyeball norm', Figure~\ref{fig:schlogl:pseudospectra} suggests that the Schl\"ogl reactions \eqref{eq:tri:molecular:matrix} are less sensitive to perturbations than the other two models, which may be found in \cite{ShevCauchyIntegralMasterEqnPseudoSpectraCTAC2015,IserlesMacNamara2017}.
An interesting experiment is to vary the parameters, namely the product $k_3 B_2$, as is done in \cite[Figure 10]{SargsyanMarzouk2009}.
In experiments not shown here, we find that when $k_3 B_2$ is varied, the numerical computation of the eigenvalues becomes less reliable (for example, when $k_3=1.4 \times 10^{-3}$ and $B_2$ is unchanged).

Figure~\ref{fig:MonoMolecular:FOV} and Figure~\ref{fig:BiMolecular:FOV}  also display the \textit{numerical range} of a matrix.
That is also known as the  \textit{field of values}, which we denote by $W$, and it is defined for a matrix $\mathbb{A}$, as the set of complex numbers that come from a quadratic form with the matrix
\begin{equation}
W(\mathbb{A}) \equiv \{x^* \mathbb{A} x \in \mathbb{C} \, : \, ||x||_2 = 1 \}.
\label{def:FOV}
\end{equation}
The field of values always contains the eigenvalues. 
A simple algorithm \cite{CharlesJohnsonFOV1978} for computing $W(\mathbb{A}) $ involves repeatedly `rotating' the matrix and then finding the \textit{numerical abscissa}
\begin{equation}
\frac{1}{2}\max(\textrm{eig}\left(\mathbb{A}+\mathbb{A}^* \right)).
\label{eq:numerical:abscissa}
\end{equation}
%It is conceivable that such an algorithm might be adapted to estimate the field of values more cheaply for the purposes of contour integrals, but we do not explore such an approach here.
The $\epsilon$-pseudospectrum of a matrix is contained in an $\epsilon$-neighbourhood of the field of values, in a sense that can be made a precise theorem \cite{TreEmb05}.
We see an example of this in Figure~\ref{fig:MonoMolecular:FOV} and Figure~\ref{fig:BiMolecular:FOV}.

\section{A Mittag-Leffler Stochastic Simulation Algorithm}
\label{sec:ML:ssa}
In this short section we provide a method for Monte Carlo simulation of the solutions of the stochastic processes that we describe. 
This Monte Carlo method can be considered  an alternative to contour integral methods, that we describe later.
As the Monte Carlo methods do not seem to fail in the presence of non-normality, they can be a useful cross-check on contour integral methods.

Here is a \texttt{MATLAB} listing to simulate Monte Carlo sample paths of the  monomolecular model of isomerization corresponding to the matrix in 
\eqref{eq:mono:molecular:matrix}, with Mittag-Leffler waiting times:

\begin{shaded}
\begin{verbatim}
t_final = 100;  m=10; c1=1; c2=1; v = [-1,  1;  1, -1]; 
initial_state = [m,0]';  alpha = 0.9; al1 = 1/alpha;  
alpi=alpha*pi; sinalpi=sin(alpi); cosalpi=cos(alpi); 
t = 0; x = initial_state; T = [0]; X = [initial_state];
while (t  < t_final)
    a(1) = c1*x(1); a(2) = c2*x(2); asum = sum(a); 
    r1=rand; r2=rand;  z=sinalpi/tan(alpi*r2)-cosalpi;
    tau = -(z/asum)^(al1)*log(r1);
    r = rand*asum;  j = find(r<cumsum(a),1,'first');
    x = x + v(:,j); t =t+tau; T = [T t];  X = [X x];
end
if (t > t_final)
    T(end) = t_final;  X(:,end) = X(:,end-1);
end
figure(1); stairs(T, X(1,:), 'LineWidth', 2); 
xlabel('Time'); ylabel('molecules of S_1'); 
title({'Isomerisation: a  monomolecular model';...
'Mittag-Leffler  SSA'; ['\alpha == ', num2str(alpha)]})
\end{verbatim}
\end{shaded}

\begin{figure}[t]
%\sidecaption[t]
\centering
\includegraphics[scale=.7]{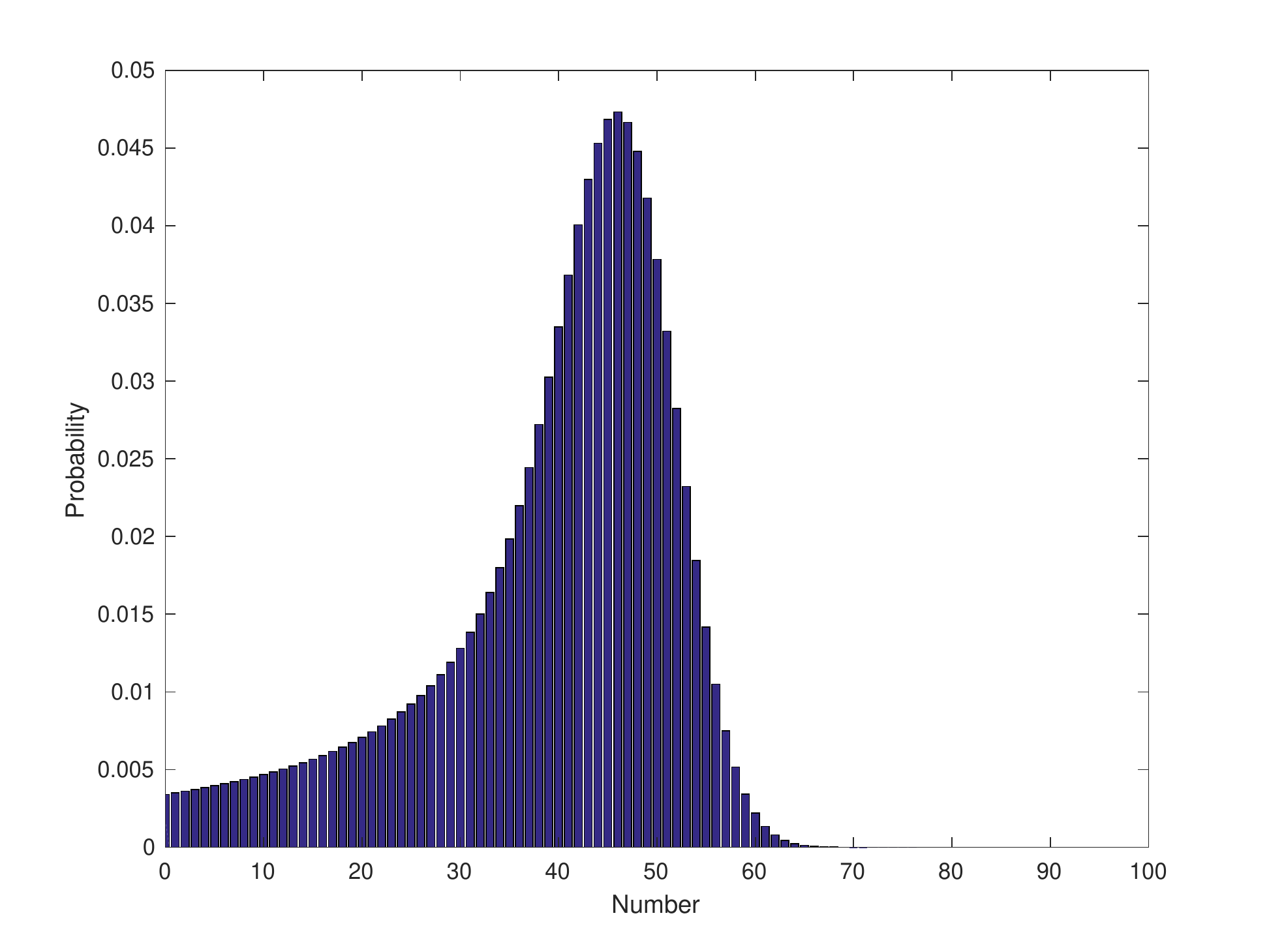}
\caption{A  discrete probability mass function is the solution of a Mittag-Leffler master equation \eqref{eq:mittag:leffler:master:equation} for the monomolecular reaction, and can be computed  with the \texttt{MATLAB} listing provided here that uses a hyperbolic contour.
Parameters: $\alpha = 0.7$, $t=1$, and initial condition a Dirac delta distribution on the state $[S_1,S_2]=[m,0]=[100,0]$. 
The $x-$axis shows the number of molecules of $S_2$.
Compare with Figure~\ref{fig:ML:SSA}.}
\label{fig:ML:ME}       % Give a unique label
\end{figure}

This is a Markov process with exponential waiting times when $\alpha=1$, in which case the program reduces to a version of the Gillespie Stochastic Simulation Algorithm.
Figure~\ref{fig:ML:SSA} shows a sample path simulated with this \texttt{MATLAB} listing.
A histogram of many such Monte Carlo samples is one way to approximate the solution of \eqref{eq:mittag:leffler:master:equation}.
A different way to compute that solution is via a contour integral, as we describe in the next section, and as displayed in Figure~\ref{fig:ML:ME}.

A time-fractional diffusion equation in the form of \eqref{eq:mittag:leffler:master:equation} in two space dimensions is solved by such a contour integral method in \cite[Figure 16.2]{TrefethenWeidemanSIAMReview2014}.
As an aside, we can modify the code above to offer a method for Monte Carlo simulation of sample paths of a closely related stochastic process, as follows.
Simplify to one space dimension (it is also straight forward to simulate in two dimensions), and by supposing that the governing $(m+1) \times (m+1)$ matrix is 
  \begin{equation}
\mathbb{A} =
\left(
\begin{tabular}{rrrrrr}
   $-1$    & $1$ &    \\
   $1$   & $-2$ &    $1$ \\
    & & $\ldots$ \\
        &    &    &    $1$ &   $-2$  &   $1$ \\
        &     &     &     &  $1$  &  $-1$
   \end{tabular}
\right)
\label{eq:one:dimensional:laplacian}
\end{equation}
with initial vector $p(0)$ being zero everywhere except the middle entry, i.e. the $\textrm{\texttt{round}}(m/2)$th entry is one. 
Then \eqref{eq:mittag:leffler:master:equation}  corresponds to  a random walk on the integers $0,1,\ldots,m$, beginning in the middle.
To simulate, modify the first few lines  of the above code segment to
\begin{verbatim}
t_final = 1;  m=10;  v = [-1,  1]; 
initial_state = round(m/2);  ...
\end{verbatim}
and also the first line in the while loop, to 
\begin{verbatim}
a(1) = (x>0); a(2) = (x<m); 
\end{verbatim}
leaving the rest unchanged.

\section{Computing a Mittag-Leffler matrix function}
\label{sec:ML:contour:integral}
In this section we first establish the utility of contour integral methods by directly computing the desired solution of  \eqref{eq:mittag:leffler:master:equation}.
This is motivation for exploring bounds on the pseudospectrum, so that  informed choices can be made when designing contour methods.

\subsection{Computing contour integrals}
\begin{figure}[t]
%\sidecaption[t]
\centering
\includegraphics[scale=.7]{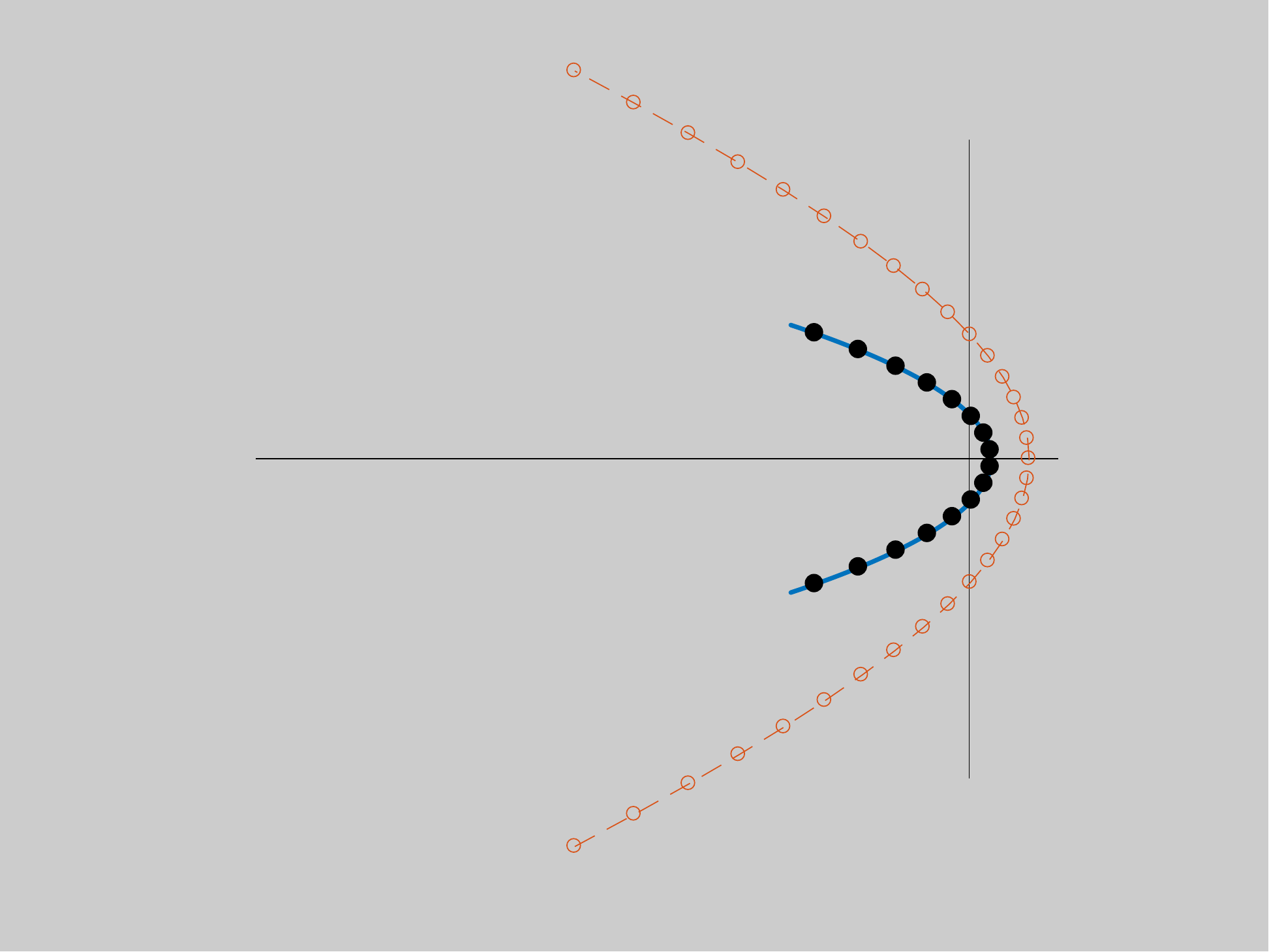}
\caption{A parabolic contour (solid) as described in \cite{TrefethenWeidemanSIAMReview2014}, and a hyperbolic contour (dashed, with $M=16$) used in the \texttt{MATLAB} listing provided here.
Nodes for a quadrature rule as in \eqref{eq:quadrature} are also marked on the contours.}
\label{fig:HankelContour}       % Give a unique label
\end{figure}

Start with \eqref{eq:mittag:leffler:master:equation}.
Take the Laplace transform. 
Then take the inverse Laplace transform.
 Of course those two steps arrive at the same solution we started with.
The advantage is  the desired solution of \eqref{eq:mittag:leffler:master:equation} represented as a contour integral \cite{ShevCauchyIntegralMasterEqnPseudoSpectraCTAC2015,MacNamaraFractionalEulerLimit2016,TrefethenWeidemanSIAMReview2014}:
%\begin{eqnarray}
%\bm{p}(t) &=& E_{\alpha}(\mathbb{A} t^{\alpha}) \bm{p}(0) = \frac{1}{2 \pi i} \int_{\Gamma} \exp(z) \left(\frac{z}{t}\mathbb{I} - \left(\frac{z}{t}\right)^{1-\alpha}\mathbb{A}\right)^{-1} \frac{\bm{p}(0)}{t} \ud z \nonumber \\
 %&=& \frac{1}{2 \pi i} \int_{\Gamma} \exp(zt) \left(z\mathbb{I} - z^{1-\alpha}\mathbb{A}\right)^{-1} \bm{p}(0) \ud z.
%\label{eq:contour:integral}
%\end{eqnarray}
\begin{equation}
\bm{p}(t) = E_{\alpha}(\mathbb{A} t^{\alpha}) \bm{p}(0)  
 = \frac{1}{2 \pi i} \int_{\Gamma} \exp(zt) \left(z\mathbb{I} - z^{1-\alpha}\mathbb{A}\right)^{-1} \bm{p}(0) \ud z.
\label{eq:contour:integral}
\end{equation}
The eigenvalues of $\mathbb{A}$ must lie to the left of the contour $\Gamma$.
In all our examples, they lie along the negative real axis.
Also, we exploit symmetry when $\mathbb{A}$ is real.
Apart from those requirements, we are free to choose the contour.
Typical choices are displayed in Figure~\ref{fig:HankelContour}.
The $\exp(zt)$ factor is nearly zero when the real part of $z$ is sufficiently negative.
This motivates choosing $\Gamma$ to go into the left-half plane because then we can neglect the infinite part of the contour that lies to the left of say, $-30$, and still obtain high accuracy.
To evaluate the integral on the remaining finite part of the contour, we use quadrature.
The trapezoidal rule is a good choice.
For a quadrature rule with $M$ nodes $z_k$ on the contour, we approximate \eqref{eq:contour:integral} by
\begin{equation}
\bm{p}(t)   \approx \sum_{k=1}^{M}  w_k(t) \bm{u}_k  
\label{eq:quadrature}
\end{equation}
where we solve a linear system for the vector $\bm{u_k}$ at each node $k$
\[
\left(z\mathbb{I} - z^{1-\alpha}\mathbb{A}\right)\bm{u_k} =  \bm{p}(0).
\]
The values $w_k(t)$ (and note that we allow these numbers to depend on $t$)  incorporate both the usual weights coming from the quadrature rule on a line, and also any scalings coming from parameterizing the contour.
There are variations of this quadrature scheme.
When $\alpha=1$, this procedure simplifies to compute the familiar matrix exponential solutions of the usual Markov systems.

%\label{sec:ML:contour:integral}
Here is a \texttt{MATLAB} listing, adapted from  Le,  Mclean \&  Lamichhane \cite{LeMcLeanLamichhane17Hyperbola}, to compute the Mittag-Leffler solution by applying the quadrature recipe  \eqref{eq:quadrature} on a \textit{hyperbolic} contour.
We choose $M=16$ quadrature points.
An example of the hyperbolic contour being used here is displayed in Figure~\ref{fig:HankelContour}.
Figure~\ref{fig:ML:ME} shows a solution of \eqref{eq:mittag:leffler:master:equation} computed with this \texttt{MATLAB} listing.

\begin{shaded}
\begin{verbatim}
function p = hyperbola_MittagLeffler(t,A,v,alpha,M)
alpha1 = 1-alpha; n = length(v); Dxi = 1.08179214/M;
xi =[-M:M]*Dxi; delta=1.17210423; mu =4.49207528*M/t;
z = mu * ( 1 - sin(complex(delta,-xi)) );
dz = 1i * mu * cos(complex(delta,-xi)) ;
c = Dxi * dz .* exp(z*t) / (2*pi*1i);
I = speye(n); p = zeros(size(v));
for k = 1:M 
    p = p + c(k)*((z(k)*I-z(k)^(alpha1)*A)\v); 
end
p = 2*real(p); k = M+1;
p = p + real(c(k)*((z(k)*I-z(k)^(alpha1)*A)\v));
\end{verbatim}
\end{shaded}

\subsection{Estimating the field of values}

One strategy to choose the contour, $\Gamma$, is to first compute a
psedospectrum of the matrix  with \texttt{EigTool} \cite{EigTool2002} (as in 
Figure~\ref{fig:schlogl:pseudospectra} for example). Then  choose  $\Gamma$ so 
that  $\|(z\mathbb{I}-\mathbb{A})^{-1}\|$ is never too large for~$z\in \Gamma$.
For example, we might choose the contour so that the 
bound~$\|(z\mathbb{I}-\mathbb{A})^{-1}\|<10^{3}$ holds (which corresponds to 
$\epsilon=10^{-3}$ in \eqref{eq:pseudospectra:definition}). Arguably, the 
value~$10^3$ could be replaced by something smaller, say $\mathcal{O}(1)$.
The particular value would depend on the application. 
In numerical experiments  with these models, the issue seems to matter only when 
$\|(z\mathbb{I}-\mathbb{A})^{-1}\|$ is significantly larger than, say, $100$.
Such an unfavourable situation can certainly arise. It is demonstrated to happen 
on the parabolic contour displayed here in Figure~\ref{fig:HankelContour} for  
bimolecular reactions~\cite{ShevCauchyIntegralMasterEqnPseudoSpectraCTAC2015}. 
The issue is also addressed by In'T Hout and 
Weideman~\cite{WeidemanBlackScholes2011} for Black--Scholes models.
Figure~\ref{fig:HankelContour} compares a parabolic contour used for a 
bimolecular model~\cite{ShevCauchyIntegralMasterEqnPseudoSpectraCTAC2015} with
a hyperbolic contour used here. A similar comparison of these contours and 
discussion can be found in a survey article by Trefethen and 
Weideman~\cite[Figure~15.2]{TrefethenWeidemanSIAMReview2014}.

However, a drawback of this procedure is that it is expensive to first compute 
the pseudospectrum. It would therefore be preferable to instead find an estimate 
of the region where the resolvent is large by some more efficient means.
We find one such estimate next.

Trefethen and Embree~\cite{TreEmb05} discuss various ways to bound the 
pseudospectrum. One approach uses the fact that if $z$ is not inside the field 
of values of a matrix, then the norm of the resolvent is bounded by the distance 
to the field of values:
\[
\|(z\mathbb{I}-\mathbb{A})^{-1}\|_2
    <\frac{1}{\mathrm{dist}\bigl(z,W(\mathbb{A})\bigr)}.
\]
This result suggests an idea for adapting the contour: choose the contour to be 
outside of the field of values. 

We now estimate the field of values of the spatial operator in the 
Fokker--Planck PDE~\eqref{eq:fokker-planck}; a similar approach has been 
applied to the Black--Scholes equation of 
finance~\cite[Theorem~3.1]{WeidemanBlackScholes2011}. We focus on the particular 
example of the  monomolecular model, but we believe this approach will be 
extended to the other models in future work.

Begin by writing the Fokker--Planck equation in the form of a conservation law,
\[
u_t+Au=0\quad\text{for $0<x<m$ and $t>0$,}
\]
where, using a dash for $\partial/\partial x$,
\begin{equation}\label{eq:second:notation}
Au=-\tfrac12(bu)''+(au)'=-\bigl(\tfrac12 bu'+Bu\bigr)'
\quad\text{and}\quad B(x)=-a(x)+\tfrac12 b'(x).
\end{equation}
For the monomolecular model, the coefficients are
\[
b(x)=c_1x+c_2(m-x)\quad\text{and}\quad a(x)=-c_1x+c_2(m-x).
\]
Here, $c_1$, $c_2$~and $m$ are positive constants.
Note that $A$ is uniformly elliptic because $b(x)\ge \min (c_1,m c_2)$ for $0 < x < m$. 

We impose either homogeneous Dirichlet boundary conditions,
\begin{equation}\label{def:dirichlet}
u(0)=0=u(m),
\end{equation}
or else zero-flux boundary conditions,
\begin{equation} \label{def:zero:flux:neumann}
\frac{1}{2}bu' + Bu=0\quad\text{for $x\in\{0,m\}$.}
\end{equation}
The domain of~$A$ is then the complex vector space~$D(A)$ of $C^2$~functions
$v:[0,m]\to\mathbb{C}$ satisfying the chosen boundary conditions.

Denote the numerical range (field of values) of~$A$ by
\begin{equation}\label{def:FOV:continuous}
W(A)=\{\,\iprod{Au,u}:\text{$u\in D(A)$ with $\iprod{u,u}=1$}\,\},
\end{equation}
where $\iprod{u,v}=\int_0^m u\bar v$. Compare this definition of~$W(A)$ in the 
continuous setting with the definition~\eqref{def:FOV} in the discrete setting.
For either \eqref{def:dirichlet}~or \eqref{def:zero:flux:neumann}, integration 
by parts gives
\[
\iprod{Au,u}=\frac12\int_0^m b|u'|^2+\int_0^mBu \bar u'.
\]
We write
\[
X+iY\equiv\iprod{Au,u}=\tfrac12P-Q
\quad\text{where}\quad
\]
\[
P=\int_0^mb|u'|^2\quad\text{and}\quad Q=-\int_0^mBu\bar{u'},
\]
and assume that 
\[
\frac{B(x)^2}{b(x)}\le K\quad\text{and}\quad
0<\beta_0\le\frac{B'(x)}{2}\le\beta_1\quad\text{for $0\le x\le m$.}
\]
In the case of zero-flux boundary conditions~\eqref{def:zero:flux:neumann}, we 
require the additional assumption
\[
B(0)\le 0\le B(m).
\]
Then we claim that for Dirichlet boundary conditions~\eqref{def:dirichlet},
\begin{equation}
\label{def:zero:dirichlet:bound}
Y^2 \le 2K(X+\beta_1)-\beta_0^2,
\end{equation}
whereas for zero-flux boundary conditions~\eqref{def:zero:flux:neumann}, 
\begin{equation}
\label{def:zero:flux:neumann:bound}
Y^2 \le 2K(X+\beta_1).
\end{equation}

To derive these estimates, first observe that
\begin{eqnarray}
2\Re Q &=& Q+\bar Q=-\int_0^m B(u\bar{u'}+\bar{u}u') \nonumber \\
	&=& -\int_0^m B(u\bar u)'=-\bigl[B|u|^2\bigr]_0^m+\int_0^m B'|u|^2.
\end{eqnarray}
Assume that $\langle u,u \rangle =1$.
The bounds on~$B'$ give
\[
\beta_0-\bigl[\tfrac12B|u|^2\bigr]_0^m 
    \le\Re Q\le\beta_1-\bigl[\tfrac12B|u|^2\bigr]_0^m,
\]
and by the Cauchy--Schwarz inequality,
\begin{eqnarray}
|Q|^2\le\biggl(\int_0^m B^2|u'|^2\biggr)\biggl(\int_0^m|u|^2\biggr)
	&=& \int_0^m\frac{B^2}{b}\,b|u'|^2  \nonumber \\
	  & \le & \biggl(\max_{[0,m]}\frac{B^2}{b}\biggr)
		\int_0^m b|u'|^2 \le KP; \nonumber
\end{eqnarray}
thus,
\[
Y^2 =(-\Im Q)^2 = |Q|^2  -(\Re Q)^2 \le KP -(\Re Q)^2.
\]
For Dirichlet boundary conditions we have $\bigl[B|u|^2\bigr]_0^m=0$, so 
$\beta_0\le\Re Q\le\beta_1$ and hence
\[
P = 2X+2\Re Q\le2X+2\beta_1\quad\text{and}\quad Y^2\le KP-\beta_0^2,
\]
implying \eqref{def:zero:dirichlet:bound}. 
For zero-flux boundary conditions, 
the extra assumptions on $B$ ensure that $\bigl[B|u|^2\bigr]_0^m \ge 0$, so 
$\Re Q \le \beta_1$ and hence
\[
P = 2X+2\Re Q\le2X+2\beta_1\quad\text{and}\quad Y^2 \le KP,
\]
implying~\eqref{def:zero:flux:neumann:bound}.

In the simple case $c_1=c_2=c$ we have 
\[
b(x)=cm\quad\text{and}\quad B(x)= -a(x) = c(2x-m)
\]
and therefore
\[
\frac{B(x)^2}{b(x)}=\frac{c}{m}(2x-m)^2\le cm\quad\text{and}\quad B'(x)= 2c
\quad\text{for $0\le x\le m$,}
\]
giving $K=cm$~and $\beta_0=\beta_1=c$. In addition, $B(0)=-cm\le0$~and 
$B(m)=cm\ge0$. Thus, for zero-flux boundary conditions,
\[
Y^2 \le 2cm(X+c).
\]

\textbf{Note that the sign convention} in the notation in~\eqref{eq:second:notation} makes the operator~$A$ positive definite (though not symmetric), whereas our  model matrices are negative definite, so $W(A)$ provides an approximation  for~$W(-\mathbb{A})=-W(\mathbb{A})$. 
We therefore have to flip signs in  \eqref{def:zero:dirichlet:bound}~and \eqref{def:zero:flux:neumann:bound} to get estimates for~$W(\mathbb{A})$, as in the following bound.

~\\

\begin{shaded}
\textbf{Parabolic estimate} of the field of values for the   \textbf{monomolecular model} matrix when $c_1=c_2=c$  in~\eqref{eq:mono:molecular:matrix}:
\begin{equation}
Y^2 \le  2cm(c-X).
\label{eq:mono:molecular:matrix:parabolic:bound}
\end{equation}
\end{shaded}

\begin{figure}[t]
%\sidecaption[t]
\centering
\includegraphics[scale=.8]{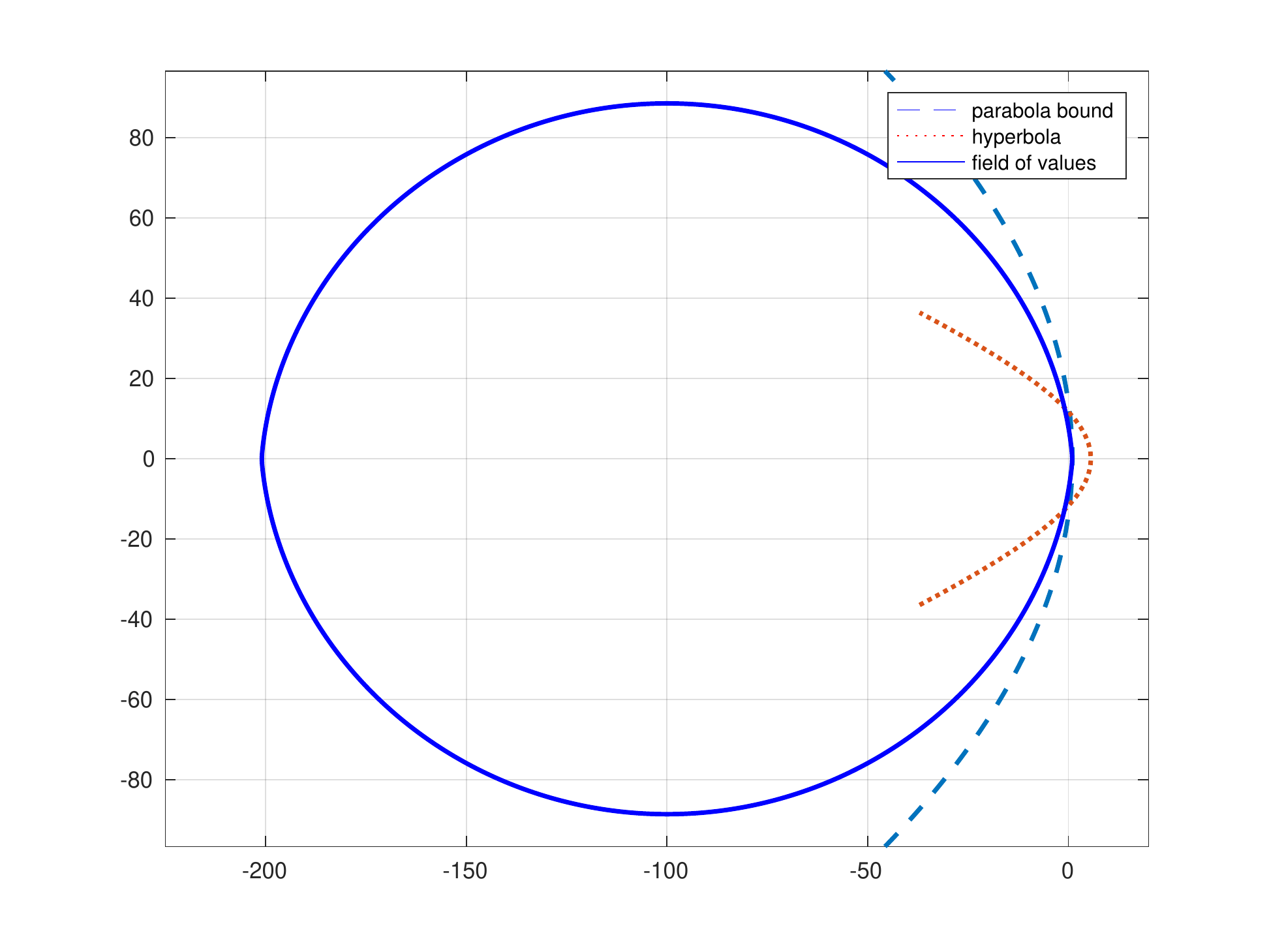}
\caption{The field of values for the discrete matrix in \eqref{eq:mono:molecular:matrix} when $c_1=c_2=1$ and $m=100$, as computed by \texttt{Chebfun}.
The dashed line is a parabolic bound in \eqref{eq:mono:molecular:matrix:parabolic:bound} for the field of values of the corresponding continuous PDE \eqref{eq:fokker-planck}.
This parabolic bound is semi-analytic and requires some analysis of the equation `by hand.'
%The figure shows numerical evidence that the latter can usefully be employed as a bound on the former.
Also displayed is a hyperbolic contour, as in Figure~\ref{fig:HankelContour} as used in the \texttt{MATLAB} listing provided here.
}
\label{fig:parabolic:bound}       % Give a unique label
\end{figure}

This bound is displayed in Figure~\ref{fig:parabolic:bound}.
As in the Figure, our matrix examples typically have a unique zero eigenvalue, with all other eigenvalues having negative real part, and the numerical abscissa  \eqref{eq:numerical:abscissa} is typically strictly positive. 
Having derived a bound in the continuous setting, we can only regard it as an approximation in the discrete setting.
Nonetheless, in these numerical experiments it does indeed seem to offer a useful bound for the discrete case.

\textbf{Discussion.}
A natural next step is to incorporate this bound into the design of the contour for methods such as the \texttt{MATLAB} listing provided in Section~\ref{sec:ML:contour:integral}.
That will be pursued elsewhere.
This article has thus laid the foundations for such attractive future research.

It is worth commenting that if a matrix is real symmetric (unlike the examples in this article) then the pseudospectrum is not an issue so it would not be good to make the contour wider (and thus also longer), and instead previously proposed (and often optimized and shorter) contours such as surveyed in~\cite{TrefethenWeidemanSIAMReview2014} would be good choices.
Making the contour wider and thus longer does impact the numerical method, but this is unavoidable in applications where the behaviour of the pseudospectra is an issue, such as the examples discussed here.
It is also worth commenting on numerical methods for computing the field of values.
Here we used \texttt{Chebfun}, which in turn is based on Johnson's algorithm~\cite{CharlesJohnsonFOV1978}.
It is conceivable that a different algorithm might be devised to estimate the field of values more cheaply for the purposes of contour integrals, but we do not explore such an approach here.
However, we do observe in these numerical experiments, as displayed in the figures for example, that the field of values --- and therefore also the estimate that we derive from it --- seems to be too conservative.
That is especially noticeable in regions where the real part is large and negative.
It might therefore also be worth considering methods of more cheaply, directly, estimating the pseudospectra.
We began with the continuous PDE so using a coarse mesh in a PDE-solver, such as a finite difference method, might be one way to obtain a cheap estimate of the pseudospectra. 

%It remains to be seen if that observation is significant for contour integral methods, for which such regions are typically  negligible due to the exponential factor $e^z$ (that approaches zero fast as the real part of $z$ becomes increasingly negative) in the integrand.

\section{Conclusion}
We have discussed, and illustrated by example, the significance of the pseudospectra for stochastic models in biology, including both Markovian and non-Markovian generalisations.
Although we focused exclusively on the $2$-norm,  the pseudospectra of these stochastic processes are perhaps more naturally addressed in the $1$-norm.
In that regard, future work will explore ways to incorporate the $1$-norm methods described by Higham and Tisseur \cite{HighamTissuerMatrixOneNorm}, to adaptively choosing the contour.

Numerical methods to compute the matrix exponential functions and matrix Mittag-Leffler functions via contour integrals must take the pseudospectrum of the matrix into account.
In particular, such methods must choose contours that are wide enough, and that adapt to avoid regions of the complex plane where the norm of the resolvent matrix is too large.
We have derived a simple, parabolic bound on the field of values of an associated Fokker--Planck PDE, which can be used as an approximation to the field of values of the corresponding discrete matrix model.
We believe this bound can help inform us when making a good choice for the contour.
Ultimately, how to devise contours in a  truly adaptive fashion, so that they can be efficiently computed, automatically and without a priori analysis, remains an important open question.

\paragraph{Acknowledgement}
The authors thank the organisers of the Computational Inverse Problems theme at the MATRIX, 2017.

\clearpage
\bibliographystyle{siam}
\bibliography{refs}
\end{document}